\documentclass[12pt]{article}

\usepackage{amssymb}
\usepackage{amsfonts}
\usepackage{amsmath}
\usepackage{amsthm}
\usepackage{graphics}
\usepackage{latexsym}

\newtheorem{theorem}{Theorem}[section]
\newtheorem{lemma}[theorem]{Lemma}
\newtheorem{remark}[theorem]{Remark}

\newtheorem{proposition}[theorem]{Proposition}
\newtheorem{definition}[theorem]{Definition}
\newtheorem{corollary}[theorem]{Corollary}

\newcommand{\bZ}{\mathbb{Z}}
\newcommand{\bN}{\mathbb{N}}
\newcommand{\bR}{\mathbb{R}}

\newcommand{\Var}{\mbox{\rm Var\,}}
\newcommand{\Cov}{\mbox{\rm Cov\,}}

\title{A central limit theorem for the sample autocorrelations of a L\'evy driven continuous
time moving average process}
\date{}

\author{Serge Cohen\thanks{Institut de Math\'ematiques de Toulouse, Universit\'e Paul Sabatier, Universit\'e de Toulouse, 118 route de Narbonne F-31062 Toulouse Cedex  9. E-mail: \texttt{Serge.Cohen@math.univ-toulouse.fr}} \and Alexander Lindner\thanks{Institut f\"ur Mathematische Stochastik,
Technische Universit\"at Braunschweig, Pockelsstra{\ss}e 14, D-38106
Braunschweig, Germany. E-mail: \texttt{a.lindner@tu-bs.de}
(Corresponding author)}}

\begin{document}
\maketitle

\abstract{In this article we consider L\'evy driven continuous time
moving average processes observed on a lattice, which are stationary
time series. We show asymptotic normality of the sample mean, the
sample autocovariances and the sample autocorrelations. A comparison
with the classical setting of discrete moving average time series
shows that in the last case a correction term should be added to the
classical Bartlett formula that yields the asymptotic variance. An
application to the asymptotic normality of the estimator of the
Hurst exponent of fractional L\'evy processes is also deduced from
these results.}
\bigskip


\section{Introduction}

 Statistical models are often written in a continuous time setting for theoretical
reasons (e.g. diffusions). But if one wants to estimate the
parameters of these models, one usually assumes only the observation
of a discrete sample. At this point a very general question, the
answer of which depends on the model chosen, is to know if the
estimation should not have been performed with an underlying
discrete model in the beginning. In this article we will consider
this for moving average processes and we refer to the classical
moving average time series models as a discrete counterpart of this
continuous model.

To be more specific, let $L = (L_t)_{t\in \bR}$ be a two sided one-dimensional L\'evy
process, i.e. a stochastic process with independent and stationary
increments, c\`adl\`ag sample paths and which satisfies $L_0=0$.
Assume further that $L$ has finite variance and expectation zero,
and let $f:\bR\to \bR$ be in $L^2(\bR)$. Let $\mu\in \bR$. Then the
process $(X_t)_{t\in \bR}$, given by
\begin{equation} \label{eq-ma1}
X_t  =  \mu + \int_{\bR} f(t-s) \, dL_s , \quad t\in \bR,
\end{equation}
can be defined in the $L^2$ sense and is called a {\it continuous
time moving average process with mean $\mu$ and kernel function $f$,
driven by $L$.}  See also~\cite{Cohen12} for more information on
such processes, in particular fractional L\'evy processes.
The process $(X_t)_{t\in \bR}$ is
then strictly stationary.
 Equation~\eqref{eq-ma1} is the natural continuous time
analogue of discrete time moving average processes
\begin{equation} \label{eq-ma3}
\widetilde{X}_t = \mu + \sum_{i\in \bZ} \psi_{t-i} Z_i, \quad t\in
\bZ,
\end{equation}
where $(Z_t)_{t\in \bZ}$ is an independent and identically
distributed (i.i.d.) noise sequence with finite variance and
expectation zero, and $(\psi_i)_{i\in\bZ}$ is a square summable
sequence of real coefficients. The asymptotic behaviour of the
sample mean and sample autocorrelation function of $\widetilde{X}_t$
in \eqref{eq-ma3} has been studied for various cases of noise
sequences $(Z_i)_{i\in \bZ}$, such as regularly varying noise
(cf.~Davis and Mikosch~\cite{Davis98}), martingale difference
sequences (cf.
Hannan~\cite{Hannan76}), 
or i.i.d. sequences with finite fourth moment or finite variance but
more restrictive conditions on the decay of the sequence
$(\psi_i)_{i\in \bZ}$ (cf. Section~7 of Brockwell and
Davis~\cite{Brockwell87}).

Another approach to obtain limit theorems for
sample autocovariances is to prove strong mixing properties of the
time series under consideration, and provided it has finite
$(4+\delta)$-moment, use the corresponding central limit theorems
(such as in Ibragimov and Linnik~\cite{Ibragimov71}, Theorem~18.5.3). If even
stronger strong mixing conditions hold, then existence of a fourth
moment may be enough. Observe however that processes with long
memory are often not strongly mixing, and in this paper we are
aiming also at applications with respect to the fractional L\'evy noise,
which is not strongly mixing.

In this paper we shall study the asymptotic behaviour as
$n\to\infty$ of the sample mean \begin{equation}
\label{eq-sample-mean} \overline{X}_{n;\Delta} := n^{-1}
\sum_{i=1}^n X_{i\Delta},
\end{equation}
of the process $(X_t)_{t\in \bR}$ defined in \eqref{eq-ma1} when
sampled at $(\Delta n)_{n\in \bN}$, where $\Delta > 0$ is fixed,
and of its sample autocovariance and sample autocorrelation function
\begin{eqnarray}
\widehat{\gamma}_{n;\Delta} (\Delta h) & := & n^{-1}
\sum_{i=1}^{n-h} (X_{i\Delta} - \overline{X}_{n;\Delta})
(X_{(i+h)\Delta} -\overline{X}_{n;\Delta}), \quad h\in \{0, \ldots,
n-1\},\quad \quad
\label{eq-sample-cov1}\\
\widehat{\rho}_{n;\Delta} (\Delta h) & := &
\widehat{\gamma}_{n;\Delta}(\Delta h) /
\widehat{\gamma}_{n;\Delta}(0), \quad h \in \{0,\ldots, n-1\}.
\label{eq-sample-corr1}
\end{eqnarray}
We write $\bN=\{0,1,2,\ldots\}$.  Under appropriate conditions on
$f$ and $L$, in particular assuming $L$ to have finite fourth moment
for the sample autocorrelation functions, it will be shown that
$\overline{X}_{n;\Delta}$ and $(\widehat{\rho}_{n;\Delta} (\Delta),
\ldots, \widehat{\rho}_{n;\Delta}(h\Delta))$ are asymptotically
normal for each $h\in \bN$ as $n\to\infty$. This is similar to the
case of discrete time moving average processes of the form
\eqref{eq-ma3} with i.i.d. noise, but unlike for those, the
asymptotic variance of the sample autocorrelations of model
\eqref{eq-ma1} will turn out to be given by Bartlett's formula {\it
plus} an extra term which depends explicitly on the fourth moment of
$L$, and in general this extra term does not vanish. This also shows
that the ``naive'' approach of trying to write the sampled process
$(X_{n\Delta})_{n\in \bZ}$ as a discrete time moving average process
as in \eqref{eq-ma3} with i.i.d. noise does not work in general,
since for such processes the asymptotic variance would be given by
Bartlett's formula only. If $\mu=0$, then further natural estimators
of the  autocovariance and autocorrelation are given by
\begin{eqnarray}
{\gamma}_{n;\Delta}^* (\Delta h) & := & n^{-1} \sum_{i=1}^{n}
X_{i\Delta}  X_{(i+h)\Delta} , \quad h\in \{0, \ldots, n-1\},
\label{eq-sample-cov2}\\
{\rho}_{n;\Delta}^* (\Delta h) & := & {\gamma}_{n;\Delta}^*(\Delta
h) / {\gamma}_{n;\Delta}^*(0), \quad h \in \{0,\ldots, n-1\},
\label{eq-sample-corr2}
\end{eqnarray}
and the conditions we have to impose to get asymptotic normality of
$\gamma_{n;\Delta}^*$ and $\rho_{n;\Delta}^*$ are less restrictive
than those for $\widehat{\gamma}_{n;\Delta}$ and
$\widehat{\rho}_{n;\Delta}$.

 We will be particularly interested in the
case when $f$ decays like a polynomial, which is e.g. the case for
fractional L\'evy noises. For a given L\'evy process with expectation
zero and finite variance, and a parameter $d\in (0,1/2)$, the {\it
(moving average) fractional L\'evy process $(M_{t;d})_{t\in \bR}$
with Hurst parameter $H:= d + 1/2$} is given by
\begin{equation}
\label{fLp} M^1_{t;d} := \frac{1}{\Gamma (d+1)}
\int_{-\infty}^\infty \left[ (t-s)^d_+ - (-s)^d_+ \right] \,
dL_s,\quad t\in \bR\end{equation} (cf.
Marquardt~\cite{Marquardt06}).  A process also called fractional
L\'evy process was introduced before by Benassi et
al.~\cite{Cohen03}, where $(x)_+= \max(x,0) $ is replaced by an
absolute value in~\eqref{fLp},
\begin{equation}
\label{fLp2} M^2_{t;d}:=\int_{-\infty}^\infty \left[
|t-s|^d - |s|^d \right]
 \, dL_s,\quad t\in \bR.\end{equation}
Although both processes have different distributions, they enjoy similar
properties. For instance the sample paths of both versions
are H\"older continuous, have the same pointwise H\"older exponent, and they
are both locally self-similar (see~\cite{Cohen03} for the definition
of this local property of their distributions). The corresponding {\it fractional
L\'evy noises} based on increments of length $\Delta > 0$ are given by
$$X^i_{t} = M^i_{t; d} - M^i_{t-\Delta; d}, \quad t\in \bR \quad i =1,\;2.$$
Hence the fractional L\'evy noise is a L\'evy driven moving average
process with kernel function \begin{equation} f^1_{d,\Delta} (s) = \frac{1}{\Gamma (d+1)}
\left ( s_+^d - (s-\Delta)_+^d \right ), \quad s \in \bR, \label{eq-kernal-fl}
\end{equation}
or
\begin{equation} f^2_{d,\Delta} (s) =
|s|^d - |s-\Delta|^d, \quad s \in \bR. \label{eq-kernal-f2}
\end{equation}
While the kernel functions $f^i_{d,\Delta}$, $ i\in \{1, 2\}$, do
not satisfy the assumptions we will impose for the theorems
regarding the sample mean $\overline{X}_{n;\Delta}$ and the sample
autocorrelation function $\widehat{\rho}_{n;\Delta}$, for $d\in
(0,1/4)$ they do satisfy the assumptions we impose for the
asymptotic behaviour of $\rho_{n;\Delta}^*$, so that an
asymptotically normal estimator of the autocorrelation and hence of
the Hurst index can be obtained if $d\in (0,1/4)$. For general $d\in
(0,1/2)$, one may take the differenced fractional L\'evy noises
$M^i_{t;d} - 2 M^i_{t-\Delta; d} + M^i_{t - 2 \Delta ;d}$, $t\in
\bR$, and our theorems give asymptotically normal estimators for the
autocorrelation function of these processes. Please note that
asymptotically normal estimators of the Hurst exponent for $
M^2_{t;d}$ are already described in~\cite{Cohen03} but they use
fill-in observations of the sample paths $ X^2(k/2^n)$ for $k = 1,
\dots, 2^n-1.$ If  $ L $ is a Brownian motion, then   $ X^1
\stackrel{d}{=} C X^2 ,$ where $\stackrel{d}{=}$ means equality in
distribution for processes and $ C $ is a constant, is the
fractional Brownian motion and it is self-similar. Except in this
case, fractional L\'evy processes are not self-similar and therefore
observations on a grid $ k/2^n $ do not yield the same information
as the time series $ X^i(t)$, $t \in \bZ.$

The paper is organised as follows: in the next section we will
derive asymptotic normality of the sample mean. Then, in
Section~\ref{S-3} we will derive central limit theorems for the
sample autocovariance  $\widehat{\gamma}_{n;\Delta}$ and the sample
autocorrelation $\widehat{\rho}_{n;\Delta}$, as well as for the
related estimators $\gamma^*_{n;\Delta}$ and $\rho^*_{n;\Delta}$ of
\eqref{eq-sample-cov2} and \eqref{eq-sample-corr2}. As a byproduct
of the asymptotic normality, these quantities are consistent
estimators of the autocovariance and autocorrelation. 
Section~\ref{S-6} presents an application of our
 results to the estimation of the parameters of fractional L\'evy
noises, where the underlying Hurst parameter is estimated. We also
recall there that fractional L\'evy noises are mixing in the
ergodic-theoretic sense, and we prove that
they fail to be strongly mixing. 

Throughout the paper, unless indicated otherwise, $L$ will be a
L\'evy process with mean zero and finite variance $\sigma^2 = E
L_1^2$, and $X=(X_t)_{t\in \bR}$ denotes the process defined in
\eqref{eq-ma1} with kernel $f\in L^2(\bR)$, $f:\bR\to\bR$. Its
autocovariance at lag $h\in \bR$ will be denoted by
\begin{equation}\label{eq-gamma2}  \gamma (h) = \gamma_f(h) = \Cov (X_0, X_{h}) = \sigma^2
\int_{\bR} f(-s) f(h-s) \, ds,
\end{equation}
where the last equation follows from the It\^o isometry.

Let us set some notations used in the sequel.

If $ v$ is vector or $ A $ a matrix  the transposed is denoted by $
v'$, respectively by $ A'.$

Convergence in distribution is denoted by $\stackrel{d}{\to}$.

The function  $\mathbf{1}_A$ for a set $A$ is one for $ x \in A,$
and vanishing elsewhere.

The autocorrelation of $X$ at lag $h$ will be denoted by $\rho (h) =
\rho_f(h) = \gamma(h)/\gamma(0).$

\section{Asymptotic normality of the sample mean} \label{S-2}
\setcounter{equation}{0}

The sample mean $\overline{X}_{n;\Delta}$ of the moving average
process $X$ of \eqref{eq-ma1} behaves like the sample mean of a
discrete time moving average process with i.i.d. noise, in the sense
that it is asymptotically normal with variance $\sigma^2
\sum_{k=-\infty}^\infty \gamma (k\Delta)$, provided the latter is
absolutely summable.

\begin{theorem} \label{thm-sample-mean}
Let $L$ have zero mean and variance $\sigma^2$, let $\mu\in \bR$ and
$\Delta > 0$. Suppose that \begin{equation}
\label{eq-absolute-summability} \left( F_\Delta : [0,\Delta]\to
[0,\infty], \quad u \mapsto F_\Delta (u) = \sum_{j=-\infty}^\infty
|f(u+j\Delta)|\right) \in L^2([0,\Delta]).
\end{equation} Then $\sum_{j=-\infty}^\infty |\gamma (\Delta j)|
<\infty$,
\begin{equation} \label{eq-gamma1}
\sum_{j=-\infty}^\infty \gamma(\Delta j) = \sigma^2 \int_0^\Delta
\left( \sum_{j=-\infty}^\infty f(u+j \Delta) \right)^2 \, du,
\end{equation}
and the sample mean of $X_\Delta,\ldots, X_{n\Delta}$ is
asymptotically normal as $n\to\infty$, more precisely
$$\sqrt{n} \, \overline{X}_{n;\Delta} \stackrel{d}{\to} N\left(\mu,
\sigma^2 \int_0^\Delta \left( \sum_{j=-\infty}^\infty f(u+\Delta j)
\right)^2  \, du\right) \quad \mbox{\rm as $n\to\infty$.}$$
\end{theorem}
\begin{remark}
Throughout the paper the assumption that $ L $ has  zero mean can be
dropped very often. For instance, if  $f\in L^1(\bR)\cap L^2(\bR)$,
then the assumption of zero mean of $L$ presents no restriction, for
in that case $L'_t := L_t - t E(L_1)$, $t\in\bR$, defines another
L\'evy process with mean zero and the same variance, and it holds
\begin{equation} \label{eq-ma2}
X_t = \mu + E(L_1) \int_{\bR} f(s) \, ds + \int_{\bR} f(t-s) \,
dL_s', \quad t\in \bZ,
\end{equation}
which has mean $\mu + E(L_1) \int_{\bR} f(s) \, ds$.
\end{remark}

\begin{proof}
For simplicity in notation, assume that $\Delta = 1$, and write $F =
F_1$. Continue $F$ periodically on $\bR$ by setting
$$F(u) = \sum_{j=-\infty}^\infty |f(u+j)| ,\quad u \in \bR.$$
Since
$$|\gamma_f(h)| \leq \sigma^2 \int_{-\infty}^\infty |f(-s)| \,
|f(h-s)| \, ds$$ by \eqref{eq-gamma2}, we have
\begin{eqnarray}
\label{uppbound-cov}
\frac{1}{\sigma^2} \sum_{h=-\infty}^\infty |\gamma_f(h)| & \leq &
\int_{-\infty}^\infty |f(-s)| \sum_{h=-\infty}^\infty |f(h-s)| \, ds
\nonumber
\\
& = & \int_{-\infty}^\infty |f(s)| F(s) \, ds \nonumber \\
& = & \sum_{j=-\infty}^\infty \int_0^1 |f(s+j)| F(s) \, ds \nonumber \\
& = & \int_0^1 F(s) \, F(s) \, ds < \infty.
\end{eqnarray}
The same calculation without the modulus gives \eqref{eq-gamma1}.

The proof for asymptotic normality is now much in the same spirit as
for discrete time moving average processes, by reducing the problem
to $m$-dependent sequences first and then applying an appropriate
variant of Slutsky's theorem. By subtracting the mean we may assume
without  loss of generality  that $\mu=0$. For $m\in \bN$, let $f_m := f \,
\mathbf{1}_{(-m,m)}$, and denote
$$X_t^{(m)} :=
\int_{\bR} f_m(s)\, dL_s = \int_{t-m}^{t+m} f(t-s)\, dL_s ,   \quad
t\in \bZ.$$
 Observe that $(X_{t}^{(m)})_{t\in\bZ}$ is a
$(2m-1)$-dependent sequence, i.e. $(X_{j}^{(m)})_{j\leq t}$ and
$(X_{j}^{(m)})_{j\geq t+2m}$ are independent for each $t\in\bZ$.
From the central limit theorem for strictly stationary
$(2m-1)$-dependent sequences (cf. Theorem~6.4.2 in Brockwell and
Davis~\cite{Brockwell87}) we then obtain that
\begin{equation} \label{eq-slutsky1}
\sqrt{n}\, \overline{X}_{n;1}^{(m)} = n^{-1/2} \sum_{t=1}^n
X_{t}^{(m)} \stackrel{d}{\to} Y^{(m)}, \quad n\to\infty,
\end{equation}
where $Y^{(m)}$ is a random variable such that
$$Y^{(m)} \stackrel{d}{=} N(0,v_m)$$
with $v_m = \sum_{j=-2m}^{2m} \gamma_{f_m}(j)$. Since
$\lim_{m\to\infty} \gamma_{f_m} (j) = \gamma_f(j)$ for each $j\in
\bZ$ by \eqref{eq-gamma2}, since
$$|\gamma_{f  \mathbf{1}_{(-m,m)}}(j)| {\leq} \sigma^2 \int_{-\infty}^\infty  |f(-s)|\, |f(j-s)| \, ds,$$  and
$\sum_{j=-\infty}^\infty \int_{-\infty}^\infty  |f(-s)|\, |f(j-s)| < \infty $ by~\eqref{uppbound-cov},
it follows from
Lebesgue's dominated convergence theorem that $\lim_{m\to\infty} v_m
= \sum_{j=-\infty}^\infty \gamma_f (j)$. Hence by \eqref{eq-gamma1},
\begin{equation} \label{eq-slutsky2}
Y^{(m)} \stackrel{d}{\to} Y, \quad m\to\infty, \quad \mbox{where}
\quad Y \stackrel{d}{=} N\left(0, \sigma^2 \int_0^1 \left(
\sum_{j=-\infty}^\infty f(u+j)\right)^2 \, du\right).
\end{equation}
A similar argument gives $\lim_{m\to\infty} \sum_{j=-\infty}^\infty
\gamma_{f-f_m} (j) = 0$, so that
\begin{eqnarray*}
\lefteqn{\lim_{m\to\infty} \lim_{n\to\infty} \Var \left( n^{1/2}
(\overline{X}_{n;1} - \overline{X}_{n;1}^{(m)}) \right)}\\
 & = &
\lim_{m\to\infty} \lim_{n\to\infty} n \, \Var \left( n^{-1}
\sum_{t=1}^n \int_{-\infty}^\infty (f(t-s)- f_m(t-s)) \,dL_s
\right),\\
& = & \lim_{m\to\infty} \sum_{j=-\infty}^\infty \gamma_{f-f_m} (j) =
0,
\end{eqnarray*}
where we used Theorem~7.1.1 in Brockwell and
Davis~\cite{Brockwell87} for the second equality.
 An
application of Chebychef's inequality then shows that
\begin{equation*} \label{eq-slutsky3} \lim_{m\to\infty}
\limsup_{n\to\infty} P( n^{1/2} |\overline{X}_{n;1} -
\overline{X}_{n;1}^{(m)}| > \varepsilon) = 0
\end{equation*} for every $\varepsilon > 0$. Together with
\eqref{eq-slutsky1} and \eqref{eq-slutsky2} this implies the claim
by a variant of Slutsky's theorem (cf.~\cite{Brockwell87},
Proposition~6.3.9).
\end{proof}
\begin{remark}\label{rem-FL1}
Let us start with an easy remark on a necessary condition on
the kernel $ f $ to apply the previous theorem.
Obviously $  F_{\Delta} \in L^1([0,\Delta]) $ is equivalent to
$ f \in L^1(\mathbb R).$ Hence
$ F_{\Delta} \in L^2([0,\Delta]) \Rightarrow f \in L^1(\mathbb R).$
\end{remark}
\begin{remark} \label{rem-expl-condition1}
Unlike for the discrete time moving average process of
\eqref{eq-ma3}, where absolute summability of the autocovariance
function is guaranteed by absolute summability of the coefficient
sequence, for the continuous time series model \eqref{eq-ma1} it is
not enough to assume that the kernel satisfies $f\in L^1(\bR) \cap
L^2(\bR)$. An example is given by taking $\Delta=1$ and
$$f(u) := \begin{cases} 0, & u \leq 0, \\
1, & u \in [0,1),\\
\frac{1\cdot 3 \cdots (2j-1)}{2^j j!} (u-j)^j ,&  u \in [j,j+1),
\quad j\in \bN.\end{cases}$$ For then the function $F_1$ is given by
$$F_1(u) = \sum_{j\in \bZ} f(u+j) = (1-u)^{-1/2}, \quad u \in
[0,1),$$ so that $F_1 \in L^1([0,1]) \setminus L^2([0,1])$. But
$F_1\in L^1([0,1])$ is equivalent to $f\in L^1(\bR)$, and since
$|f(u)|\leq 1$ for all $u\in \bR$, this implies also $f\in
L^2(\bR)$. Observe further that for non-negative $f$, condition
\eqref{eq-absolute-summability} is indeed necessary and sufficient
for absolute summability of the autocovariance function.
\end{remark}

\section{Asymptotic normality of the sample autocovariance}
\label{S-3} \setcounter{equation}{0}

 As usual, we consider the stationary
process \begin{equation} \label{eq-def-X} X_t =
\int_{-\infty}^\infty f(t-s) \, dL_s, \quad t\in \bR.\end{equation} 
We recall that
$$\gamma^*_{n;\Delta} (h\Delta) = n^{-1} \sum_{t=1}^n X_{t\Delta} X_{(t+h)\Delta}, \quad h \in
\bN $$ and first we establish an asymptotic result for $\Cov
(\gamma^*_{n;\Delta}(p \Delta ), \gamma^*_{n;\Delta}(q \Delta)).$


\begin{proposition} \label{prop-cov1}
Let $L$ be a (non-zero) L\'evy process, with expectation zero, and
finite fourth moment, and denote $\sigma^2 := E L_1^2$ and $\eta :=
\sigma^{-4} E L_1^4$.  Let $\Delta > 0,$ and suppose further that $f
\in L^2(\bR) \cap L^4(\bR)$ and that \begin{equation}
\label{eq-cyprus} \left( [0,\Delta] \to \bR, \quad u \mapsto
\sum_{k=-\infty}^\infty f(u+k \Delta)^2 \right) \in L^2([0,\Delta]).
\end{equation}
For $q\in \bZ$ denote
$$g_{q;\Delta} : [0,\Delta] \to \bR, \quad u \mapsto \sum_{k=-\infty}^\infty
f(u+k \Delta) f(u+(k+q)\Delta ), $$ which belongs to  $ L^2([0,\Delta]),$
by the previous assumption. If further
\begin{equation} \label{eq-cov-summe}
\sum_{h=-\infty}^\infty |\gamma(h \Delta)|^2 < \infty,
\end{equation}
then we have  for each $p,q \in \bN$
\begin{multline}
 \label{eq-cov-limit}
 \lim_{n\to\infty} n \Cov (\gamma^*_{n;\Delta}(p \Delta ), \gamma^*_{n;\Delta}(q \Delta))
 =  (\eta -3)\sigma^4 \int_0^\Delta g_{p;\Delta}(u) g_{q;\Delta}(u) \, du +\\
\sum_{k=-\infty}^\infty \big[ \gamma(k\Delta ) \gamma((k-p+q)\Delta ) + \gamma((k+q)\Delta)
\gamma((k-p)\Delta) \big].
\end{multline}
\end{proposition}

\begin{proof}
For simplicity in notation we assume that $\Delta = 1$. The general
case can be proved analogously or reduced to the case $\Delta=1$ by
a simple time change. We shall first show that for $t,p,h,q\in \bZ$
\begin{eqnarray}
\lefteqn{E (X_t X_{t+p} X_{t+h+q} X_{t+h+p+q})} \nonumber\\
& = & (\eta-3)\sigma^4 \int_{-\infty}^\infty f(u) f(u+p) f(u+h+p)
f(u+h+p+q) \, du \nonumber\\
& & + \gamma(p) \gamma(q) + \gamma(h+p) \gamma(h+q) + \gamma(h+p+q)
\gamma(h). \label{eq-cov-darst1}
\end{eqnarray}
To show this, assume first that $f$ is of the form
\begin{equation}
\label{eq-def-special-f} f(s) = f_{m,\epsilon} (s) =
\sum_{i=-m/\epsilon}^{m/\epsilon} \psi_i \mathbf{1}_{(i \epsilon,
(i+1)\epsilon]}(s), \end{equation} where $m\in\bN$, $\epsilon>0$
such that $1/\epsilon \in \bN$, and $\psi_i\in \bR$,
$i=-m/\epsilon,\ldots, m/\epsilon$. Denote
$$X_{t; m,\epsilon} := \int_{-\infty}^\infty f_{ m,\epsilon} (t-s) \, dL_s =
\sum_{i=-m/\epsilon}^{m/\epsilon} \psi_i (L_{t-i\epsilon} -
L_{t-(i+1)\epsilon}),\quad t\in \bR.$$ Denote further
$$Z_i := L_{i \epsilon} - L_{(i-1) \epsilon}, \quad i \in\bZ.$$
Then $(Z_i)_{i\in\bZ}$ is i.i.d. and we have
$$X_{t\epsilon; m , \epsilon} = \sum_{i=-m/\epsilon}^{m/\epsilon} \psi_i
Z_{t-i}, \quad t\in \bZ.$$

At this point we will need to compute the fourth moment
of integrals of the L\'evy process. Let us state an elementary result
that yields a formula for this moment.
\begin{lemma}
\label{lem:4-moment} Let  $\phi\in L^2(\bR) \cap L^4(\bR)$, then,
with the assumptions and notations on $L$ used in
Proposition~\ref{prop-cov1},
\begin{equation}
  \label{eq:mom}
  E (\int_{\bR} \phi(s) dL_s)^4 = (\eta - 3)\sigma^4\ \int_{\bR} \phi^4(s) ds + 3 \sigma^4 (\int_{\bR} \phi^2(s) ds)^2.
\end{equation}
\end{lemma}
\begin{proof}
If $\nu$ is the L\'evy measure of $L$ and $A$ its Gaussian variance,
then by the L\'evy Khintchine formula we get
\begin{eqnarray*}
\xi(u) & = & E \exp \left (i u \int_{\bR} \phi(s) dL_s \right
)\\
&  = & \exp \left ( -\frac12 A u^2 \int_{\bR} \phi^2(s) ds + \int_{
\bR \times \bR} [ e^{i u \phi(s) x} - 1 - i u \phi(s) x ] \nu(dx) ds
\right ).\end{eqnarray*} Then $ E (\int_{\bR} \phi(s) dL_s)^4$ is
obtained as the fourth derivative of $\xi$ at $ u =0.$ If we recall
that $ (\eta - 3) \sigma^4 = \int_{\bR} x^4 \nu(dx),$ and $
\sigma^2= A+\int_{\bR} x^2 \nu(dx),$ we get~\eqref{eq:mom}, after
elementary but tedious computations.
\end{proof}
To continue with the proof of Proposition~\ref{prop-cov1},  we now
apply~\eqref{eq:mom} to the special case where $
f(s)=\mathbf{1}_{(0,\epsilon]}(s) $ and we get
\begin{equation} \label{eq-fourth-moment} E
Z_i^2 = E L_\epsilon^2 = \sigma^2\epsilon,\quad E Z_i^4 = E L_\epsilon^4 =
\eta\sigma^4\epsilon - 3 \sigma^4\epsilon + 3 \sigma^4 \epsilon^2.
\end{equation}
As shown in the proof of Proposition~7.3.1 in \cite{Brockwell87}, we then
have
\begin{eqnarray*}
\lefteqn{ E( X_{t; m,\epsilon} X_{t+p; m,\epsilon} X_{t+h+p; m,\epsilon}
X_{t+h+p+q; m,\epsilon})} \\
& = & \left( E Z_i^4 - 3 (E Z_i^2)^2\right)
\sum_{i=-m/\epsilon}^{m/\epsilon} \psi_i \psi_{i+p/\epsilon} \psi_{i
+ h/\epsilon + p/ \epsilon} \psi_{i+h/\epsilon +
p/\epsilon + q/\epsilon}\\
& & + \gamma_{m,\epsilon}(p) \gamma_{m,\epsilon}(q) + \gamma_{m,\epsilon}
(h+p) \gamma_{m,\epsilon} (h+q) + \gamma_{m,\epsilon} (h+p+q)
\gamma_{m,\epsilon} (h),
\end{eqnarray*}
where $\gamma_{m,\epsilon} (u) = E (X_{0; m,\epsilon} X_{u;
m,\epsilon})$, $u\in\bR$. By \eqref{eq-fourth-moment}, $$E Z_i^4 - 3
(E Z_i^2)^2 = (\eta-3)\sigma^4 \epsilon,$$ and \begin{eqnarray*}
\lefteqn{\epsilon \sum_{i=-m/\epsilon}^{m/\epsilon} \psi_i
\psi_{i+p/\epsilon} \psi_{i + h/\epsilon + p/ \epsilon}
\psi_{i+h/\epsilon + p/\epsilon + q/\epsilon}}\\ & = &
\int_{-\infty}^\infty f(u) f(u+p) f(u+h+p) f(u+h+p+q) \, du,
\end{eqnarray*}
so that \eqref{eq-cov-darst1} follows for $f$ of the form
$f=f_{m,\epsilon}$. Now let $f\in L^2(\bR) \cap L^4(\bR)$ and $X_t$,
$t\in\bR$, defined by \eqref{eq-def-X}. Then there is a sequence of
functions $(f_{m_k,\epsilon_k})_{k\in\bN}$ of the form
\eqref{eq-def-special-f} such that $f_{m_k,\epsilon_k}$ converges to
$f$ both in $L^2(\bR)$ and in $L^4(\bR)$ as $k\to\infty$.
 Then for each fixed $t\in\bR$, we
have that $X_{t; m_k,\epsilon_k} \to X_t$ in $L^2(P)$ ($P$ the
underlying probability measure) as $k\to\infty$, where we used the
It\^o isometry. Further, by Lemma~\ref{lem:4-moment}, and convergence
of $f_{m_k,\epsilon_k}$ both in $L^2(\bR)$ and in $L^4(\bR),$
we get convergence of $ X_{t; m_k,\epsilon_k}$ to $ X_t$ in  $L^4(P).$
This then shows \eqref{eq-cov-darst1}, by letting
$f_{m_k,\epsilon_k}$ converge to $f$ both in $L^2(\bR)$ and $L^4(\bR)$
and observing that $\gamma_{m_k,\epsilon_k} (u) \to \gamma (u)$ for
each $u\in\bR$. From \eqref{eq-cov-darst1} we  conclude 
that, with $p,q\in \bN$,
\begin{equation} \label{eq-Tk}
\Cov (\gamma_{n; 1}^* (p ), \gamma_{n; 1}^*(q)) = n^{-1}
\sum_{|k|<n} (1 - n^{-1} |k|) T_k,
\end{equation}
where
\begin{eqnarray*}
T_k & = & \gamma(k  ) \gamma(k-p+q) + \gamma(k+q) \gamma(k-p) \\
& & + (\eta-3) \sigma^4 \int_{-\infty}^\infty f(u) f(u+p) f(u+k)
f(u+q+k) \, du.
\end{eqnarray*}
Now by \eqref{eq-cov-summe}, $\sum_{k=-\infty}^\infty |T_k| <
\infty$ if \begin{equation} \label{eq-hilfsreihe1}
\sum_{k=-\infty}^\infty \left|\int_{-\infty}^\infty f(u) f(u+p)
f(u+k) f(u+q+k) \, du\right|\end{equation} is finite. Denote
$$G_r(u) := \sum_{k=-\infty}^\infty |f(u+k) f(u+k+r)|, \quad
u\in\bR, \quad r\in \bN.$$ Then $G_r$  is periodic, and by
assumption, $G_r$ restricted to $[0,1]$ is square integrable. Hence
we can estimate \eqref{eq-hilfsreihe1} by
\begin{eqnarray*}
& & \sum_{k=-\infty}^\infty \int_{-\infty}^\infty |f(u) f(u+p)| \,
|f(u+k) f(u+q+k) | \, du \\
& = & \sum_{h=-\infty}^\infty \int_{h}^{h+1} |f(u) f(u+p)| \, G_q(u)
\, du \\
&  =& \sum_{h=-\infty}^\infty \int_0^{1} |f(u+h) f(u+p+h)| \, G_q(u)
\, du \\
& = & \int_0^{1} G_p(u) G_q(u) \, du < \infty.
\end{eqnarray*}
The same calculation without the modulus and an application  of the
dominated convergence theorem to \eqref{eq-Tk} then shows
\eqref{eq-cov-limit}.
\end{proof}

\begin{remark}
A sufficient condition for 
\eqref{eq-cov-summe} is that
$\sum_{h=-\infty}^\infty |\gamma(h\Delta)| <\infty$, which is
implied by the function $F_\Delta$ in Theorem~\ref{thm-sample-mean}
belonging to $L^2([0,\Delta])$. Another sufficient condition is that
$\Phi: u \mapsto \sum_{k\in\bZ} |\mathcal{F} (f) (u+2\pi
k/\Delta)|^2$, $u\in [0,2\pi/\Delta]$ is in $L^\infty
([0,2\pi/\Delta])$, where $\mathcal{F} (f)$ is the Fourier transform
of $f\in L^2(\bR)$ in the form $z\mapsto \int_{-\infty}^\infty
e^{izt} f(t)\, dt$ (for $L^1$-functions). For if $\|\Phi\|_\infty
\leq B$, then $(f(\cdot +h\Delta))_{h\in\bZ}$ is a Bessel sequence
in $L^2(\bR)$ with bound $B/\Delta$, i.e.
$$\sum_{h=-\infty}^\infty \left|\int_{-\infty}^\infty \varphi(u) f(u+h\Delta)\,
du \right|^2 \leq B\Delta^{-1} \int_{-\infty}^\infty \varphi(u)^2 \,
du \quad \forall\; \varphi\in L^2(\bR),$$ see e.g. Theorem~7.2.3 in
Christensen~\cite{Christensen03}. Taking $\varphi=f$ then gives the
square summability of the autocovariance functions by
\eqref{eq-gamma2}.

Please remark that $\sum_{h=-\infty}^\infty \gamma(h\Delta)^2 <
\infty$ cannot be deduced from the condition that $u\mapsto
\sum_{k=-\infty}^\infty f(u+k\Delta)^2 $ is in $L^2([0,\Delta])$.
One can take $\Delta=1$ and $ f(s)= \sum_{i \ge 1} \frac{{\mathbf
1}_{(i,i+1]}(s)}{i^H} $ for $ \frac12 < H \leq \frac34,$ to get the
latter condition but not $ \sum_{h=-\infty}^\infty \gamma(h)^2 <
\infty.$
\end{remark}
\begin{remark} \label{rem-positive}
By \eqref{eq-gamma2}, the condition \eqref{eq-cov-summe} can be
written as \linebreak $\sum_{k=-\infty}^\infty \left(
\int_{-\infty}^\infty f(s) f(s+k\Delta) \, ds\right)^2 < \infty$.
The assumption \eqref{eq-cov-summe1} used in Theorem~\ref{th:as-n}
below is slightly stronger than \eqref{eq-cov-summe}, but equivalent
to \eqref{eq-cov-summe} if $f\geq 0$.
%
\end{remark}

The following theorem gives asymptotic normality of the sample
autocovariance and sample autocorrelation and the related estimators
$\gamma^*_{n;\Delta}$ and $\rho^*_{n;\Delta}$.

\begin{theorem}
  \label{th:as-n}
(a) Suppose  the assumptions of Proposition~\ref{prop-cov1} are
satisfied and suppose further that
\begin{equation} \label{eq-cov-summe1}
\sum_{k=-\infty}^\infty \left( \int_{-\infty}^\infty |f(s) f(s+k \Delta)|
\, ds\right)^2 < \infty.
\end{equation}
Then we have
for each $h \in \bN$
\begin{equation} \label{limit-behaviour1}
\sqrt{n} (\gamma_{n;\Delta}^* (0)-\gamma(0) , \ldots, \gamma_{n;\Delta}^*(h \Delta)-
\gamma(h))' \stackrel{d}{\to} N(0, V) , \quad n\to\infty,
\end{equation}
where $V = (v_{pq})_{p,q=0,\ldots, h} \in \bR^{h+1,h+1}$ is the
covariance matrix defined by
\begin{multline}
v_{pq}
 =  (\eta -3)\sigma^4 \int_0^\Delta g_{p;\Delta}(u) g_{q;\Delta}(u) \, du +\\
\sum_{k=-\infty}^\infty \big[ \gamma(k\Delta ) \gamma((k-p+q)\Delta ) + \gamma((k+q)\Delta)
\gamma((k-p)\Delta) \big]. \label{eq-V}
\end{multline}
(b) In addition to the assumptions of (a), assume that the function
 $$u\mapsto \sum_{j=-\infty}^\infty |f(u+j \Delta)|$$ is in
$L^2([0,\Delta])$. Denote by
$$\widehat{\gamma}_{n;\Delta}(j\Delta) = n^{-1} \sum_{t=1}^{n-j} (X_{t\Delta} - \overline{X}_{n;\Delta})
(X_{(t+j) \Delta} - \overline{X}_{n;\Delta}), \quad j=0,1,\ldots,
n-1,$$ the sample autocovariance, as defined in
\eqref{eq-sample-cov1}. Then we have for each $h\in\bN$
$$\sqrt{n} (\widehat{\gamma}_{n;\Delta} (0)-\gamma(0) ,
\ldots, \widehat{\gamma}_{n;\Delta}(h \Delta)- \gamma(h))'
\stackrel{d}{\to} N(0, V) , \quad n\to\infty,$$ where $V =
(V_{pq})_{p,q=0,\ldots, h}$ is defined by
\eqref{eq-V}.\\
(c) For $j\in\bN$ let $\rho_{n;\Delta}^*(j\Delta) =
\gamma_{n;\Delta}^*(j \Delta)/\gamma_{n;\Delta}^*(0)$ and
$\widehat{\rho}_n(j\Delta) =
\widehat{\gamma}_{n;\Delta}(j\Delta)/\widehat{\gamma}_{n;\Delta}(0)$,
the latter being the sample autocorrelation at lag $j \Delta$.
Suppose that $f$ is not almost everywhere equal to zero. Then, under
the assumptions of (a), we have for each $h\in\bN$, that
\begin{equation} \label{eq-corr-asymptotic}
\sqrt{n} (\rho^*_{n;\Delta} (\Delta) - \rho(\Delta), \ldots, \rho^*_{n;\Delta}(h \Delta) - \rho(h \Delta) )'
\stackrel{d}{\to} N(0,W),\quad n\to\infty,
\end{equation}
where $W=W_\Delta= (w_{ij;\Delta})_{i,j=1,\ldots, h}$ is given by
$$w_{ij;\Delta} = \widetilde{w}_{ij;\Delta} + \frac{(\eta-3)\sigma^4}{\gamma(0)^2}
\int_0^\Delta \big( g_{i;\Delta}(u) - \rho(i \Delta)
g_{0;\Delta}(u)\big) (g_{j;\Delta}(u) - \rho(j \Delta)
g_{0;\Delta}(u)\big) \, du,$$ and
\begin{eqnarray*} \widetilde{w}_{ij;\Delta} & = & \sum_{k=-\infty}^\infty
\big( \rho((k+i)\Delta) \rho((k+j)\Delta) + \rho((k-i)\Delta)
\rho((k+j)\Delta) + 2 \rho(i\Delta)\rho(j\Delta) \rho(k\Delta)^2 \\
& & \quad - 2\rho(i\Delta) \rho(k\Delta) \rho((k+j)\Delta) - 2\rho(j\Delta) \rho(k\Delta) \rho((k+i)\Delta)
\big) \\
& = & \sum_{k=1}^\infty \big( \rho((k+i)\Delta) + \rho((k-i)\Delta) - 2\rho(i\Delta)
\rho(k\Delta)\big) \times \\
 & & \phantom{\sum_{k=1}^\infty( \rho((k+i)\Delta) + \rho((k-i)\Delta))} \quad \big(\rho((k+j)\Delta) + \rho((k-j)\Delta) - 2\rho(j\Delta)\rho(k\Delta) \big)
\end{eqnarray*}
is given by Bartlett's formula. If additionally the function
 $u\mapsto \sum_{j=-\infty}^\infty |f(u+j\Delta)|$ is in
$L^2([0,\Delta])$, then it also holds that
\begin{equation} \label{eq-corr-asymptotic2}
\sqrt{n} (\widehat{\rho}_{n;\Delta} (\Delta) - \rho(\Delta), \ldots,
\widehat{\rho}_{n;\Delta}(h\Delta) - \rho(h\Delta) )' \stackrel{d}{\to} N(0,W),\quad
n\to\infty.
\end{equation}
\end{theorem}
\begin{proof}
For  simplicity in notation  we assume again $ \Delta=1$ in this
proof.

(a)
Using Proposition~\ref{prop-cov1} it follows as in  the proof of
Proposition~7.3.2 in \cite{Brockwell87}, that the claim is true if $f$ has
additionally compact support. 
 For general $f$ and $m\in\bN$ let $f_m
:= f  \mathbf{1}_{(-m,m)}$. Hence we have that
$$n^{1/2} ( \gamma_{n;(m)}^{*}(0) - \gamma_{m}(0) , \ldots, \gamma_{n;(m)}^*(h) -
\gamma_{m}(h))' \stackrel{d}{\to} \mathbf{Y}_m, \quad m\to\infty,$$
where $\gamma_{m}$ is the autocovariance function of the process
$X_{t;m} = \int_{-\infty}^\infty f_m(t-s) \, dL_s$,
$\gamma_{n;{(m)}}^*(p) = n^{-1} \sum_{t=1}^n X_{t;m} X_{t+p;m}$ the
corresponding autocovariance estimate, and $\mathbf{Y}_m
\stackrel{d}{=} N(0, V_{m})$ with $V_m = (v_{pq;m})_{p,q = 0,\ldots,
h}$ and $$ v_{pq;m}
 =  (\eta -3)\sigma^4 \int_0^1 g_{p;(m)}(u) g_{q;(m)}(u) \, du +
\sum_{k=-\infty}^\infty \big[ \gamma_{m}(k) \gamma_{m}(k-p+q) +
\gamma_{m}(k+q) \gamma_{m}(k-p) \big].$$ Here, $g_{p;(m)} (u) =
\sum_{k=-\infty}^\infty f_m(u+k) f_m(u+k+p)$, $u\in [0,1]$.\\
Next, we want to show that $\lim_{m\to\infty} V_m = V$. Observe
first that
$$g_{p;(m)} (u) = \sum_{k=-\infty}^\infty f_m(u+k) f_m(u+k+p) \to
\sum_{k=-\infty}^\infty f(u+k) f(u+k+p) = g_{p,1}(u)=: g_{p}(u) $$
almost surely in the variable $u$ as $m\to\infty$  by Lebesgue's
dominated convergence theorem, since $u\mapsto
\sum_{k=-\infty}^\infty |f(u+k) f(u+k+p)|$ is in $L^2([0,1])$ by
\eqref{eq-cyprus} and hence is almost surely finite. Further we have
$$|g_{p;(m)}(u)| \leq \sum_{k=-\infty}^\infty |f(u+k) f(u+k+p)|$$
uniformly in $u$ and $m$, so that again by the dominated convergence
theorem we have that $g_{p;(m)} \to g_p$ in $L^2([0,1])$ as
$m\to\infty$. Next, observe that
$$|\gamma_{m}(k)| \leq \int_{-\infty}^\infty |f(s) f(s+k)| \, ds
\quad \forall\; m\in \bN\quad \forall\; k\in\bZ.$$ Since
$\lim_{m\to\infty} \gamma_{m}(k) = \gamma(k)$ for every $k\in\bZ$,
it follows from the dominated convergence theorem and
\eqref{eq-cov-summe1} that $(\gamma_{m}(k))_{k\in\bZ}$ converges in
$l^2(\bZ)$ to $(\gamma (k))_{k\in\bZ}$. This together with the
convergence of $g_{p;(m)}$ gives the desired $\lim_{m\to\infty} V_m
= V$, so that
$$\mathbf{Y}_m \stackrel{d}{\to} \mathbf{Y}, \quad m\to\infty,$$
where $\mathbf{Y} \stackrel{d}{=} N(0,V)$. Finally, that
$$\lim_{m\to\infty} \limsup_{n\to\infty} P(n^{1/2} |\gamma_{n;(m)}^*
(p) - \gamma_{m}(p) - \gamma^*(p) + \gamma(p)| > \varepsilon) =
0\quad \forall\; \varepsilon > 0,\; p \in \{0,\ldots, h\}$$ follows
as in Equation (7.3.9) in~\cite{Brockwell87}.
An application of a variant of Slutsky's theorem
(cf.~\cite{Brockwell87}, Proposition 6.3.9) then gives the
claim.\\
(b) This follows as in  the proof of Proposition~7.3.4
in~\cite{Brockwell87}. One only has to observe that by
Theorem~\ref{thm-sample-mean}, $\sqrt{n}\, \overline{X}_{n;1}$
converges in distribution to a normal random variable as
$n\to\infty$. In particular, $\overline{X}_{n;1}$ must
converge to 0 in probability as $n\to\infty$.\\
(c) The limit theorem follows as in the proof of Theorem~7.2.1 in~\cite{Brockwell87}, and for $w_{ij}$ we have the
representation
\begin{eqnarray*}
w_{ij;\Delta} & = & (v_{ij} - \rho(i) v_{0j} - \rho(j) v_{i0} +
\rho(i)
\rho(j) v_{00}) / \gamma(0)^2 \\
& = & \widetilde{w}_{ij;\Delta} + \frac{(\eta-3)\sigma^4}{\gamma(0)^2} \times\\
 & & \quad \int_{0}^1 \big( g_i(u) g_j(u) - \rho(i) g_0(u) g_j(u) - \rho(j)
g_i(u) g_0(u) + \rho(i) \rho(j) g_0(u)^2\big) \, du,
\end{eqnarray*}
giving the claim.
\end{proof}

\begin{remark}\label{rem-moving}
It is easy to check that $w_{ij;\Delta} = \widetilde{w}_{ij;\Delta}$
if $f$ is of the form $f= \sum_{i=-\infty}^\infty \psi_i
\mathbf{1}_{(i\Delta,(i+1)\Delta]}$, in accordance with Bartlett's
formula, since then $(X_{t\Delta})_{t\in\bZ}$ has a discrete time
moving average representation with i.i.d. coefficients.
\end{remark}
\begin{remark}
Another case when $w_{ij;\Delta} = \widetilde{w}_{ij;\Delta}$ is
when $\eta=3$, which happens if and only if $L$ is Brownian motion.
However, in general we do not have $w_{ij;\Delta} =
\widetilde{w}_{ij;\Delta}$. An example is given by $f =
\mathbf{1}_{(0,1/2]} + \mathbf{1}_{(1,2]}$ and $ \Delta=1,$ in which
case $g_{1;1} = \mathbf{1}_{(0,1/2]}$ and $g_{0;1} = 2 \cdot
\mathbf{1}_{(0,1/2]} + \mathbf{1}_{(1/2,1]}$, and it is easy to see
that $g_{1;1} -\rho(1) g_{0;1}$ is not almost everywhere zero, so
that $w_{11;1} \neq \widetilde{w}_{11;1}$ if $\eta\neq 3$. The
latter example corresponds to a moving average process, which is
varying at the scale $\frac12$, but sampled at integer times.
Observe however that $w_{11;1/2} = \widetilde{w}_{11;1/2}$ by
Remark~\ref{rem-moving}. A more detailed study of such phenomena in
discrete time can be found in Niebuhr and Kreiss~\cite{Niebuhr12}.
%
\end{remark}

\begin{remark}
Recently,  sophisticated and powerful results on the normal
approximation of Poisson functionals using Malliavan calculus have
been obtained. E.g., Peccati and Taqqu \cite[Theorems 2, 3 and
5]{Peccati2008} prove a central limit theorem for double Poisson
integrals and apply this to a specific quadratic functional of a
L\'evy driven Ornstein--Uhlenbeck process, and Peccati et
al.~\cite[Section 4]{Peccati2010} obtain bounds for such limit
theorems, to name just of few of some recent publications on this
subject. It may be possible to apply the results of
\cite{Peccati2010, Peccati2008} to obtain another proof of
Theorem~\ref{th:as-n} under certain conditions such as finite 6th
moment, but we have not investigated this issue further. Note that
our  proof uses only basic knowledge of stochastic integrals and
methods from time series analysis.
\end{remark}


\section{An application to fractional L\'evy noise} \label{S-6}
\setcounter{equation}{0}
We will now apply the previous results to fractional L\'evy processes.
Recall from \eqref{fLp} and \eqref{fLp2} that these were  denoted by
$$
M^1_{t;d} := \frac{1}{\Gamma (d+1)} \int_{-\infty}^\infty \left[
(t-s)^d_+ - (-s)^d_+ \right] \, dL_s,\quad t\in \bR, \quad
\mbox{and}
$$
$$ M^2_{t;d}=\int_{-\infty}^\infty \left[
|t-s|^d - |s|^d \right]
 \, dL_s,\quad t\in \bR,$$
respectively, and the corresponding {\it fractional L\'evy noises}
based on increments of length $\Delta > 0$  by
$$X^i_{t} = M^i_{t; d} - M^i_{t-\Delta; d}, \quad t\in \bR, \quad i =1,\;2.$$
Hence the fractional L\'evy noises are  L\'evy driven moving average
processes with kernel functions $$ f^1_{d,\Delta} (s)
=\frac{1}{\Gamma (d+1)} \left ( s_+^d - (s-\Delta)_+^d \right ) ,
\quad s \in \bR,$$ and
$$ f^2_{d,\Delta} (s) =
|s|^d - |s-\Delta|^d, \quad s \in \bR, $$ respectively. Neither $
f^1_{d,\Delta}$ nor $f^2_{d,\Delta}$ are in $L^1(\bR)$, so
Theorem~\ref{thm-sample-mean} cannot be applied because of
Remark~\ref{rem-FL1}. Please note that for the same reason the
assumptions for {\em (b)} of Theorem~\ref{th:as-n} and
for~\eqref{eq-corr-asymptotic2} are not fulfilled.

For  simplicity in notation  we assume $ \Delta=1, $ and drop the
subindex $\Delta.$ Although the fractional noises $ X^i$ have
different distributions for $ i = 1 $ and $ i = 2,$ they are both
stationary with the autocovariance
\begin{equation}
  \label{eq:auto-cov-fLn}
   \mathbb E (X^i_{t+h} X^i_t)  = \gamma_{X^i}(h) =
\frac{C^i(d)\sigma^2}{2} \left( |h+1|^{2 d + 1} - 2 |h|^{2 d + 1 } + |h-1|^{2 d + 1}
\right),
\end{equation}
where $ C^i(d)$ is a normalising multiplicative constant depending on $d.$
Both processes $  X^i $ are infinitely divisible and of moving average type,
hence we know from~\cite{Cambanis95, Fuchs12}
 that $(X^i_t)_{t\in\bZ}$ is mixing in
the ergodic-theoretic sense. For fixed $h\in \bZ$,
define the function
$$F: \bR^\bZ \to \bR, \quad (x_n)_{n\in\bZ} \mapsto x_0 x_h.$$
If $T$ denotes the forward shift operator, then
$$F (T^k (X^i_t)_{t\in\bZ}) = X^i_k X^i_{k+h},$$
and from Birkhoff's ergodic theorem (e.g. Ash and Gardner~\cite{Ash75},
Theorems~3.3.6 and 3.3.10) we know that
$$\frac{1}{n} \sum_{k=1}^n X^i_k X^i_{k+h} \to E \left( F(
(X^i_t)_{t\in\bZ}) \right) = E X^i_0 X^i_h,\quad n\to\infty,$$ for $
i =1,\; 2,$ and the convergence is almost sure and in $L^1$
(\cite{Ash75}, Theorems~3.3.6 and 3.3.7). 
Hence, with $\gamma_n^* = \gamma_{n;1}^*$ as defined in
\eqref{eq-sample-cov2},
$$\lim_{n\to\infty} \gamma_n^* (h)  =
\frac{C^i(d) \sigma^2}{2} \left( |h+1|^{2 d + 1 } - 2 |h|^{2 d + 1}
+ |h-1|^{2 d + 1} \right) \quad \mbox{a.s.}$$ Since $\gamma(0) =
C^i(d) \sigma^2$ and $\gamma(1) = C^i(d) \sigma^2 (2^{2 d  } - 1),$
$\rho_n^*(1)=\frac{\gamma_n^* (1)}{\gamma_n^* (0)}$ is a strongly
consistent estimator for $2^{2 d } -1$. Hence,
\begin{equation}
\label{eq:dhat}
  \widehat{d}:=\frac{1}{2} \left( \frac{ \log (\rho^*_n(1) +1)}{\log
2}\right)
\end{equation}
 is a strongly consistent estimator for $d$. 

The question of the asymptotic normality of these estimators arises
naturally. There are many classical techniques to show the
asymptotic normality of an ergodic stationary sequence by assuming
some stronger mixing assumption. As far as we know, they do not work
in our setting. To illustrate this point, we would like to show that
fractional L\'evy noises are not strongly mixing. Let us first
recall the definition.
\begin{definition}
\label{def:strong-mix} {\rm Let $(X_n)_{n \in \bZ} $ be a stationary
sequence, and let
$$ \alpha_X(n) = \sup \{ |P(A \cap B ) - P(A) P(B) |, \; A \in \sigma( X_k,\; k \le 0),\; B \in \sigma( X_k,\; k \ge n) \}.$$
The sequence $ (X_n)_{n \in \bZ} $ is {\it strongly mixing} if $
\lim_{n \to \infty} \alpha_X(n) = 0.$}
\end{definition}
In our case we know the weak mixing property $ \lim_{n \to \infty}
|P(A \cap B ) - P(A) P(B) | = 0$ for $  A \in \sigma( X_0),\; B \in
\sigma( X_n),$ because of~\cite{Cambanis95, Fuchs12}. There are
classical central limit theorems for strongly mixing sequences,
see e.g. \cite{Merlevede2006} for an overview. The following result,
which is stated as Proposition 34 in \cite{Merlevede2006}, will be
particularly useful for us.


\begin{theorem}
Suppose that $(X_t)_{t\in \bZ}$ is a mean zero, strongly mixing
sequence and that there exists some $\delta > 0$ and a constant
$K>0$ such that
\begin{eqnarray}
\label{eq-1} E |X_0|^{2+\delta} & < & \infty,\\
\label{eq-2} \lim_{m\to\infty} \Var (\sum_{i=1}^m X_i) & = & \infty, \\
\label{eq-3} E |\sum_{i=1}^m X_i|^{2+\delta} & \leq & K \left(\Var
(\sum_{i=1}^m X_i)\right)^{1+\delta/2}\quad \forall\; m\in \bN.
\end{eqnarray}
Write $$S^{(n)} (t):= \sum_{i=1}^{\lfloor nt \rfloor} X_i,\quad n\in
\bN, \quad t\in [0,1],$$
where $ \lfloor x \rfloor $ is the integer part of the real number $x.$ Then
\begin{equation} \label{eq-4}
\left( \Var(\sum_{i=1}^n X_i ) \right)^{-1/2} S^{(n)} \stackrel{d}{\to}  B
\quad \mbox{weakly in $D[0,1]$,}
\end{equation}
where $B$ is a standard Brownian motion.
\end{theorem}

\begin{lemma} \label{lem-4.3}
Let $L$ be a two sided non-zero L\'evy process with expectation zero
and finite fourth moment, and let $d\in (0,1/2)$. Let
$$X^i_t := M^i_{t;d} - M^i_{t-1;d}, \quad t\in \bZ,$$
be  the corresponding fractional L\'evy noises, for $i=1, \; 2.$
Then $(X^i_t)_{t\in \bZ}$ satisfies \eqref{eq-1} -- \eqref{eq-3} with $\delta=2$, for
$i=1, \;2.$
\end{lemma}

\begin{proof}
The proof is only written for the fractional noise $X^1$ denoted by $ X $ but it is similar
for $X^2.$
Equation~\eqref{eq-1} holds since $L$ has finite fourth moment and
since the kernel function of fractional noise is in $L^2(\bR) \cap
L^4(\bR)$, and \eqref{eq-2} follows from the fact that
\begin{equation}
  \label{eq:sum}
 \sum_{i=1}^m X_i = M_{m;d},
  \end{equation}
and
\begin{equation} \label{eq-5}
 \Var (M_{m;d}) = C m^{2d+1}\quad \forall\;
m\in \bN \end{equation} for some constant $C$. To see \eqref{eq-3} for $ \delta=2,$ we use~\eqref{eq:sum}
and
$$f_m(s) := \frac{1}{\Gamma(d+1)} [(m-s)_+^d - (-s)_+^d], \quad s\in
\bR.$$ Then by Lemma~\ref{lem:4-moment},
\begin{eqnarray*}  E |M_d(m)|^4  & =  & (\eta - 3)\sigma^4\ \int_{\bR} f_m^4(s) ds + 3 \sigma^4 \left (\int_{\bR} f_m^2(s) ds \right )^2.
\end{eqnarray*}
Observe that
$$\left( \int_\bR f_m^2(s) \, ds\right)^2 = C m^{4d+2}$$ and that
$$\int_\bR f_m^4 (s)\, ds \leq  C'  m^{4d+1} $$
for positive constants $C,\;C'$, which gives the claim.
\end{proof}

A consequence of the previous theorem and lemma is the following
negative result.

 \begin{corollary} Assume the assumptions of Lemma \ref{lem-4.3}. Then the fractional L\'evy  noises
 $X^1$ and $X^2$ are not strongly mixing.
\end{corollary}
\begin{proof}
If fractional L\'evy noises were strongly mixing, then \eqref{eq-4}
would follow. Please remark that, since fractional noises are
increments,  $ S^{(n)}(t)= M^i_d(\lfloor nt \rfloor) - M^i_d(0),$
and $ \left( \Var(\sum_{i=1}^n X_i ) \right) = C n^{2 d + 1} $
by~\eqref{eq-5}. Owing to the asymptotic self-similarity of
fractional L\'evy processes (Proposition 3.1 in~\cite{Cohen03}),  we
know that
$$ \frac{M^2_d(\lfloor nt \rfloor) - M^2_d(0)}{n^{d + 1/2}} \stackrel{d}{\to} B_{d+1/2}(t),$$
where the limit is a fractional Brownian motion with Hurst exponent
$ d+1/2.$ A similar asymptotic self-similarity holds for $ M^1_d.$
Hence \eqref{eq-4} is violated and the corollary is proved by
contradiction.
\end{proof}

 Nevertheless one can apply Theorem~\ref{th:as-n} to get the asymptotic normality
of the estimator~\eqref{eq:dhat}.
\begin{proposition}
\label{prop:AN-estd<1/4} Assume the assumptions of
Lemma~\ref{lem-4.3}, and let $\widehat{d}$ be defined by
\eqref{eq:dhat}. If $ d \in (0,1/4) $ then $\sqrt{n} (\widehat{d}
-d) $ converges in distribution to a Gaussian random variable as
$n\to\infty$.
\end{proposition}
\begin{proof}
The proof is only written for the fractional noise $X^1$ denoted by $ X $ but it is similar
for $X^2.$

We shall apply~\eqref{eq-corr-asymptotic} to get convergence of $
\sqrt{n} (\rho^*_{n} (1) - \rho(0)) $ to a Gaussian random variable.
First, observe that
\begin{equation}
  \label{eq:L2L4}
   f^1_{d,1} \in  L^2(\bR) \cap
L^4(\bR),
\end{equation}
and $\sum_{k=-\infty}^\infty \gamma(k )^2 < \infty. $ The latter
inequality is classical for the fractional Gaussian noise, when $ d
< 1/4$, and holds for fractional L\'evy noises since they have the
same autocorrelation as fractional Gaussian noise. This also implies
\eqref{eq-cov-summe1} by Remark~\ref{rem-positive} since $f_{d,1}^1
\geq 0$. Let us check that $ g_{0} := g_{0;1} \in L^2(0,1).$ Since
\begin{eqnarray*}
 \Gamma^2(d+1)  g_{0} (u)&=& \sum_{k=-\infty}^\infty ( (u + k)_+^d - (u +
 k-1)_+^d)^2,
\end{eqnarray*}
it follows for all $ u \in (0,1)$ that
\begin{eqnarray*}
  \Gamma^2(d+1) g_{0} (u)& \le & \sum_{k=0}^\infty( (1 + k )_+^d - (k-1)_+^d)^2\\
             &  =  & 1 + \sum_{k=1}^\infty |k|^{2 d} \left ( \left (1+ \frac1{|k|} \right )^d - \left (1- \frac1{|k|} \right )^d  \right )^2\\
             & <  & \infty,
\end{eqnarray*}
so that even $g_0 \in L^\infty([0,1])$. Hence the assumptions of
Theorem~\ref{th:as-n}~(a) are fulfilled, and the result follows from
\eqref{eq-corr-asymptotic}.
\end{proof}
If $ d \geq 1/4 $ then  $\sum_{k=-\infty}^\infty \gamma(k )^2 =
\infty,$ since fractional L\'evy noises have the same
autocorrelation as the fractional Gaussian noise. Hence  we consider
$$Z^i_t = X^i_t - X^i_{t-1}, \quad t\in\bZ,$$
for which it holds $\sum_{k=-\infty}^\infty \gamma_{Z^i}(k )^2 <
\infty.$
We conclude from Birkhoff's ergodic theorem that
$$\gamma_{n,Z}^* (h) := \frac{1}{n}\sum_{k=2}^n Z_k Z_{k+h} \to E (Z_0 Z_h), \quad n\to\infty,$$
is a strongly consistent estimator of
\begin{equation*}
  \label{eq:cov-Z}
  E (Z_0 Z_h)= \frac{C^i(d) \sigma^2}{2} \left( -|h+2|^{2 d + 1 }+ 4 |h+1|^{2 d + 1 }- 6  |h|^{2 d + 1} + 4 |h-1|^{2 d + 1} -|h-2|^{2 d+1}
\right).
\end{equation*}
Therefore $\rho_{n,Z}^*(1)=\frac{\gamma_{n,Z}^* (1)}{\gamma_{n,Z}^*
(0)}$ is a strongly consistent estimator for $\phi(d)=\frac{-3^{2d +
1}+ 4 \cdot 2^{2d + 1}-7}{8 - 2^{2d + 2}}.$ It turns out that $
\phi$ is increasing on $(0,1/2).$ Therefore we can define the
estimator
\begin{equation}
\label{eq:dtilde}
  \tilde{d}:= \phi^{-1}(\rho^*_{n,Z}(1)).
\end{equation}
\begin{proposition}
\label{prop:AN-estd>1/4} Assume the assumptions of
Lemma~\ref{lem-4.3}, and let $\tilde{d}$ be defined by
\eqref{eq:dtilde}.
 If $ d \in
(0,1/2) $ then $\sqrt{n} (\tilde{d} -d) $ converges in distribution
to a Gaussian random variable as $n\to\infty$.
\end{proposition}
\begin{proof}
The proof is only written for the fractional noise $X^1$ denoted by
$ X $ but it is similar for $X^2.$ Let us remark that $ Z_t=
\int_{-\infty}^\infty \tilde{f}^1_{d,1}(t-s) \, dL_s, $ where
$$   \tilde{f}^1_{d,1}(s)=  f^1_{d,1}(s)-f^1_{d,1}(s-1) , $$
To apply Theorem~\ref{th:as-n}, we have to check that
$$   \tilde{f}^1_{d,1}(s)  \in  L^2(\bR) \cap
L^4(\bR),$$ which is obvious from \eqref{eq:L2L4}. Moreover we
already know that  $\sum_{k=-\infty}^\infty \gamma_{Z^i}(k )^2 <
\infty.$ This time, however, the kernel function $\tilde{f}^1_{d,1}$
is not nonnegative, but it is easy to see that
$|\tilde{f}^1_{d,1}(t)| \leq C \min(1, |t|^{d-2})$ and hence that
$$\int_{-\infty}^\infty |\tilde{f}^1_{d,1} (t) \, \tilde{f}^1_{d,1}
(t+k) | \, dt \leq C'\min(1, |k|^{d-1}), \quad \forall\; k\in \bZ,$$
for some constants $C, C'$, giving \eqref{eq-cov-summe1}. Finally,
\begin{eqnarray*}
  \Gamma^2(d+1)  g_{0} (u)&=& \sum_{k=-\infty}^\infty ( (u + k)_+^d + (u + k-2)_+^d - 2 (u + k-1)_+^d)^2,
\end{eqnarray*}
and estimating the summands separately for $k<0$, $k=0,1$ and $k\geq
2$ we obtain for $u\in [0,1]$
\begin{eqnarray*}
  \Gamma^2(d+1) g_{0} (u)& \le &
  1 + (2^d +2)^2 + \sum_{k=2}^\infty (k^d - (k-1)^d)^2 \\
          & <  & \infty,
\end{eqnarray*}
so that $g_0 \in L^\infty([0,1]) \subset L^2([0,1])$. The claim now
follows from Theorem~\ref{th:as-n}, using
\eqref{eq-corr-asymptotic}.
\end{proof}

\subsubsection*{Acknowledgement}
Major parts of this research where carried out while AL was visiting
the Institut de Math\'ematiques de Toulouse as guest professor in
2008 and 2009. He takes pleasure in thanking them for the kind
hospitality and their generous financial support.

\def\cprime{$'$}


\end{document}